\documentclass{amsproc}
\usepackage{graphicx}
\usepackage{amscd}
\usepackage{amsmath}

\theoremstyle{plain}
\newtheorem{theorem}{Theorem}
\newtheorem{definition}{Definition}
\newtheorem{remark}{Remark}
\newtheorem{problem}{Problem}
\newtheorem{acknowledgement}{Acknowledgement}

\begin{document}
\title[On factorization of operators ]{On factorization of operators between Banach spaces}
\author{Eugene Tokarev}
\address{B.E. Ukrecolan, 33-81,Iskrinskaya str., 61005, Kharkiv-5, Ukraine }
\email{evtokarev@yahoo.com}

\begin{abstract}
A finite-dimensional analogue $GL_{fin}(X)$ of the known Gordon-Lewis
constant $GL(X)$ of a Banach space $X$ is introduced; in its definition are
used only finite rank operators. It is shown that there exist Banach spaces $%
X$ such that $GL_{fin}(X)$ is finite and $GL(X)$ is infinite. Moreover, for
any Banach space $X$ of cotype 2, the finite-dimensional constant $%
GL_{fin}(X)<\infty$.
\end{abstract}

\maketitle

\section{Introduction}

Let $X,Y$ be Banach spaces; $u:X\rightarrow Y$ be a (linear, bounded)
operator. This operator is said to be $L_{1}$-factorable if there exists a
measure $\mu$ and operators $v:X\rightarrow L_{1}(\mu)$; $%
w:L_{1}(\mu)\rightarrow Y^{\ast\ast}$ such that $k_{Y}\circ u=v\circ w$,
where $k_{Y}$ is a canonical embedding of $Y$ into $Y^{\ast\ast}$. It is
clear that for any $L_{1}$-factorable operator $u$
\begin{equation*}
\gamma_{1}(u):=\inf\{\left\| v\right\| \left\| w\right\| :k_{Y}\circ
u=v\circ w\}<\infty,
\end{equation*}
where the infimum is taken over all possible factorizations of $u$.

Grothendieck asked ([1], problem 2, p. 78):

\textit{Whether there exists absolutely summing operators which are not }$%
L_{1}$\textit{-factorable?}

Gordon and Lewis [2] answered this question in affirmative. At the same time
they showed that there exists a wide class of Banach spaces $X$ with the
following property:

(GL)\ \textit{Every absolutely summing operator }$u:X\rightarrow Y$\textit{\
where }$Y$\textit{\ is an arbitrary Banach space is }$L_{1}$\textit{%
-factorable. }

This class contains all \textit{spaces with the local unconditional structure%
} (shortly, l.u.st.), which were introduced in [2]. Recall that according to
[3] a Banach space $X$ is said to has the l.u.st. if and only if its second
conjugate is isomorphic to a complemented subspace of some Banach lattice.
Particularly, the class of all Banach spaces having the l.u.st. contains all
Banach lattices. The result of [2] further was improved (see e.g. [4], p.
107); now it is known that if a Banach space $X$ is of cotype 2 and has the
l.u.st. then each its subspace enjoy the property (GL).

Gordon and Lewis [2] defined for a given Banach space $X$ a constant
\begin{equation*}
GL\left( X\right) =\sup\{\gamma_{1}\left( u\right) /\pi_{1}\left( u\right)
:u:X\rightarrow Y\},
\end{equation*}
where the supremum is taken over all Banach spaces $Y$ and all absolutely
summing operators $u:X\rightarrow Y$ (here and below $\pi_{1}\left( u\right)
$ denotes the absolutely summing norm of an operator $u$). In their solution
of the aforementioned Grothendieck problem a sequence $\left( u_{n}\right)
_{n<\infty}$ of finite rank operators ($u_{n}:X\rightarrow l_{2}^{(m_{n})}$)
was constructed in such a way that $\sup\{\gamma_{1}\left( u_{n}\right)
/\pi_{1}\left( u_{n}\right) :n<\infty\}=\infty$.

Let us define a\textit{\ finite dimensional version of the Gordon-Lewis
constant:}
\begin{equation*}
GL_{fin}(X)=\sup\{\gamma_{1}\left( u\right) /\pi_{1}\left( u\right) :\text{ }%
u:X\rightarrow l_{2};\text{ }\dim u(X)<\infty\}.
\end{equation*}

It is clear that $GL_{fin}(X)=\infty$ implies that $GL(X)=\infty$.

However in a general case may be a situation when for a given Banach space $X
$ constants $GL_{fin}(X)$ and $GL(X)$ are essentially different. To be more
exact, it this article will be shown that for some Banach spaces $X$ the
constant $GL_{fin}(X)$ is finite and $GL(X)$ is infinite.

The question of the full description of a class \textbf{GL} of those Banach
spaces that has the property (GL) (i.e. those Banach spaces $X$ for which $%
GL(X)<\infty$) is open. In the article it would be shown that a class of all
Banach spaces $X$ for which $GL_{fin}(X)<\infty$ contains the class of all
Banach spaces of cotype 2.

\section{Definitions and notations}

Let $\mathcal{B}$ denotes a class of \textit{all} Banach spaces. For $X,Y\in%
\mathcal{B}$ let $B(X,Y)$ be a Banach space of all bounded linear operators $%
u:X\rightarrow Y$ (equipped with a strong operator topology).

Let $1\leq p\leq\infty$. An operator $u\in B(X,Y)$ is said to be

\begin{itemize}
\item  $p$-\textit{absolutely summing,} if there is a constant $\lambda>0$
such that
\begin{equation*}
(\sum\nolimits_{j<n}\left\| u\left( x_{i}\right) \right\|
^{p})^{1/p}\leq\lambda(\sum\nolimits_{j<n}\left| \left\langle
x_{i},f\right\rangle \right| ^{p})^{1/p}
\end{equation*}
for any $f\in X^{\ast}$ and any finite set $\{x_{i}:i<n;n<\infty\}\subset X$.
\end{itemize}

Its $p$-\textit{absolutely summing norm }$\pi_{p}\left( u\right) $ is the
smallest constant $\lambda$.

\begin{itemize}
\item  $p$-\textit{factorable,} if there exists a such measure $\mu$ and
operators $v\in B(X,L_{p}(\mu))$; $w\in B(L_{p}(\mu),Y^{\ast\ast})$ such
that $k_{Y}\circ u=w\circ v$, where $k_{Y}$ is the canonical embedding of $Y$
in its second conjugate $Y^{\ast\ast}$, and the symbol $\circ$ denotes the
composition of operators.
\end{itemize}

Its $\gamma_{p}$\textit{-norm }is given by
\begin{equation*}
\gamma_{p}\left( u\right) =\inf\{\left\| v\right\| \left\| w\right\|
:k_{Y}\circ u=w\circ v\}.
\end{equation*}

\begin{itemize}
\item  $p$-\textit{integral,} if there exists a such probability measure $\mu
$ and such operators $v\in B(X,L_{\infty}(\mu))$, $w\in B(L_{p}(\mu),Y^{\ast
\ast})$ that $k_{Y}\circ u=w\circ\varphi\circ v$, where $\varphi$ is an
inclusion of $L_{\infty}(\mu)$ into $L_{p}(\mu)$.
\end{itemize}

Its $p$-\textit{integral norm} is given by
\begin{equation*}
\iota_{p}(u)=\inf\{\left\| v\right\| \left\| \varphi\right\| \left\|
w\right\| :k_{Y}\circ u=w\circ\varphi\circ v\}.
\end{equation*}

\begin{itemize}
\item  \textit{Nuclear,} if there are sequences $(f_{n})_{n<\infty}\in
X^{\ast }$; $\left( y_{n}\right) _{n<\infty}\in Y$ such that for all $x\in X$
\begin{equation*}
u\left( x\right) =\sum\nolimits_{n<\infty}f_{n}\left( x\right) y_{n}.
\end{equation*}
\end{itemize}

Its \textit{nuclear norm} $\nu_{1}(u)$ is given by
\begin{equation*}
\nu_{1}\left( u\right) =\inf\{\sum\left\| f_{n}\right\| \left\|
y_{n}\right\| :u\left( x\right) =\sum\nolimits_{n<\infty}f_{n}\left(
x\right) y_{n}\},
\end{equation*}
where the infimum is taken over all possible representations of $u$.

A Banach space $X$ is said to be \textit{of cotype 2} if there exists such a
constant $c<\infty$ that for any $n<\infty$ and any subset $%
\{x_{i}:i<n\}\subset X$%
\begin{equation*}
\left( \sum_{i<n}\left\| x_{i}\right\| ^{2}\right) ^{1/2}\leq c\int
_{0}^{1}\left\| \sum_{i<n}r_{k}\left( t\right) x_{k}\right\| dt,
\end{equation*}
where $r_{k}\left( t\right) $ are Rademacher functions.

It is well known that the space $L_{p}$ is of cotype 2 for all $p\in
\lbrack1,2]$.

\section{Main results}

\begin{definition}
Let $Y$ be a Banach space of cotype 2. It will be said that $X\in\mathcal{B}$
is the \textit{Pisier's space for }$Y$, if following conditions are satisfy.

\begin{enumerate}
\item  $Y$ is isometric to a subspace of $X.X$ and $X^{\ast}$ both are of
cotype 2.

\item  Any operator $u:X\rightarrow l_{2}$ is absolutely summing; any
operator $v:X^{\ast}\rightarrow l_{2}$ is absolutely summing.

\item  There is a constant $\lambda_{X}$ such that for any finite rank
operator $u:X\rightarrow X$ its nuclear norm $\nu_{1}(u)$ satisfies the
estimate $\nu_{1}(u)\leq\lambda_{X}\left\| u\right\| $.

\item  If $X$ contains a subspace isomorphic to $l_{1}$ then its conjugate $%
X^{\ast}$ contains a subspace isomorphic to $l_{2}$.

\item  Any operator $w:l_{2}\rightarrow X$ is 2-absolutely summing; the same
is true for every operator $s:l_{2}\rightarrow X^{\ast}$.
\end{enumerate}
\end{definition}

The famous result of J. Pisier [5] shows the existence of the corresponding
Pisier's space $X$ for any Banach space $Y$ of cotype 2.

\begin{theorem}
There exists a Banach space $X$ such that $GL\left( X\right) =\infty$ and $%
GL_{fin}(X)<\infty$.
\end{theorem}

\begin{proof}
Let $X\in\mathcal{B}$ be a space which is the Pisier's space for $L_{1}[0,1]$%
. Certainly, it has subspaces that are isometric to $l_{1}$ and to $l_{2}$.
Let $l_{2}$ be isometric to a subspace $X_{2}\hookrightarrow X$ and $%
i:l_{2}\rightarrow X_{2}$ be the corresponding isometry. Recall that $X$ is
of cotype 2 and that each operator $u\in B(X,l_{2})$ is absolutely summing.
Let $u:X\rightarrow l_{2}$. The composition $w=i\circ u$ acts on $X$, is
absolutely summing and is $L_{1}$-factorable if and only if $u$ is. Assume
that $u$ has a finite dimensional range, i.e. that dim$\left( u(X)\right)
<\infty$. Then $w$ is also a finite rank operator. Hence, its nuclear norm $%
\nu\left( w\right) $ may be estimated as $\nu\left( w\right) \leq
\lambda\left\| w\right\| $, where $\lambda$ does not depend on $w$. Thus,
\begin{equation*}
\gamma_{1}\left( u\right) =\gamma_{1}\left( w\right) \leq\iota_{1}\left(
w\right) \leq\nu\left( w\right) \leq\lambda\left\| w\right\| \leq
\lambda\pi_{1}\left( w\right) .
\end{equation*}

This sequence of inequalities implies that
\begin{equation*}
\gamma_{1}\left( u\right) \leq\lambda\pi_{1}\left( u\right) .
\end{equation*}

Recall that $u$ was an arbitrary operator of finite dimensional range.
Hence, $GL_{fin}(X)<\infty$.

Let us show that for a special choosing of $u:X\rightarrow l_{2}$ such
inequalities may be not longer true. Indeed, according to [5], if $X$
contains a subspace isomorphic to $l_{1}$ then its conjugate $X^{\ast}$
contains a subspace (say, $Y_{2}$) which is isomorphic to $l_{2}$. Let $%
j:l_{2}\rightarrow Y_{2}$ be the corresponding isomorphism; $\left\|
j\right\| \left\| j^{-1}\right\| <\infty$. Because of reflexivity, $Y_{2}$
is weakly* closed and, hence, there exists a quotient map $%
h:X\rightarrow(Y_{2})^{\ast}$ which generates a surjection $H:X\rightarrow
l_{2}$: $H=(j_{2})^{\ast}\circ h$. From our assumptions on properties of $X$
follows that the conjugate operator $H^{\ast}$ is $L_{\infty}$-factorable
and is 2-absolutely summing. Clearly, $H^{\ast}$ induces an isomorphism
between $l_{2}$ and $H^{\ast}(l_{2})$. However this property contradicts
with $\pi_{2}\left( H^{\ast }\right) <\infty$. Hence, $H$ cannot be factored
through $L_{1}\left( \mu\right) $.
\end{proof}

\begin{theorem}
For every Banach space of cotype 2, $GL_{fin}(X)$ is finite.
\end{theorem}

\begin{proof}
Let a Banach space $Y$ is of cotype 2. It can be isometricaly embedded in a
corresponding Pisier's space $X_{Y}$, which has the same properties as was
listed above. Let $T:Y\rightarrow l_{2}$ be an absolutely summing operator.
Then $T$ is 2-absolutely summing and, consequently, 2-integral, i.e. admits
a factorizations $T=k\circ j\circ s$, where $s:X\rightarrow L_{\infty}(\mu)$%
; $k\circ j:L_{\infty}(\mu)\rightarrow l_{2}$. Because of injectivity of $%
L_{\infty}(\mu)$, $s$ can be extended to $\overline{s}:X_{Y}\rightarrow
L_{\infty}(\mu)$. The corresponding operator $\overline{T}=k\circ j\circ%
\overline{s}:X_{Y}\rightarrow l_{2}$ is 2-absolutely summing too and extends
$T$. Recall that from properties of $X_{Y}$ follows that $\overline {T}$ is
absolutely summing. Assume that $T$ has a finite dimensional range. The same
is true for $\overline{T}$. From theorem 1 follows that $\overline {T}$ is $%
L_{1}$-factorable. Hence $T$ is also $L_{1}$-factorable.
\end{proof}

\begin{remark}
From [6] follows that if there exists some $p>2$ such that the space $l_{p}$
is finitely representable in $X\in\mathcal{B}$ then $X$ contains a subspace $%
Z$ for which $GL_{fin}(Z)$ is infinite.

\begin{problem}
Whether every Banach space $Z$ which is not of cotype 2 contains a subspace
which does not enjoy this property?
\end{problem}
\end{remark}

\begin{acknowledgement}
Author wish to indebted his thanks to V.M. Kadets for his attention to this
paper
\end{acknowledgement}

\section{References}

\begin{enumerate}
\item  A. Grothendieck, \textit{R\'{e}sum\'{e} de la th\'{e}orie
m\'{e}trique des produits tensoriels} \textit{topologiques,} Bol. Soc. Math.
Sa\~{o} Paulo, \textbf{8} (1956), 1-79

\item  Y. Gordon and D.R. Lewis, \textit{Absolutely summing operators and
local} \textit{unconditional structures}, Acta Math., \textbf{133 }(1974),
27-48

\item  T. Figiel, W.B. Johnson and L. Tzafriri, \textit{On Banach lattices
and spaces having local unconditional structure with applications to Lorentz
function spaces}, J. Approx. Theory, \textbf{13:4} (1975), 392-412

\item  J. Pisier, \textit{Factorization of linear operators and geometry of
Banach spaces}, Publ.: Conf. Board of Math. Sci. by Amer. Math. Soc.,
Providence, Rhode Island, 1986, 153 p.

\item  J. Pisier, \textit{Counterexample to a conjecture of Grothendieck},
-Acta Math. \textbf{151} (1983) 181-208

\item  T. Figiel, S. Kwapie\'{n} and A. Pe\l czy\'{n}ski, \textit{Sharp
estimates for the constants of unconditional structure of Minkowski spaces},
Bull. Acad. Pol. Sci., S\'{e}r. Math. Astron. et Phys., \textbf{25} (1977),
1221-1226
\end{enumerate}

\end{document}